\documentclass[12pt,a4paper]{article}
\usepackage{euscript,amsfonts,amssymb,amsmath,amscd}

\input{epsf}

\newcommand{\op}{\operatorname}
\newcommand{\m}{\mathbb}
\newcommand{\e}{\iota}
\newcommand{\eps}{\varepsilon}
\begin{document}
\title{Quasitoric Manifolds with Invariant Almost Complex Structure}

\author{Andrei Kustarev}

\maketitle

\begin{abstract}

We prove that any quasitoric manifold $M$ admits a $T^n$-invariant
almost complex structure if and only if $M$ admits a positive
omniorientation. In particular, we show that all obstructions
to existence of $T^n$-invariant almost complex structure
on $M$ arise from cohomology of underlying polytope - and hence are trivial.

\end{abstract}

\section{Introduction}

Let $M^{2n}$ be an oriented closed compact manifold. We say that $M^{2n}$ is
{\it quasitoric} over simple polytope $P$
(\cite{davis},\cite{buchpanov}) if:
\begin{enumerate}
{\item there is locally standart action of torus $T^n$ on $M^{2n}$
and} {\item there exists a map $\pi:M^{2n}\to P^{n}$ which is a
superposition of projection $M^{2n}\to M^{2n}/T^n$ and
diffeomorphism of manifolds with corners $M^{2n}/T^n\to P$.}
\end{enumerate}

Any smooth projective toric variety gives the example of quasitoric manifold,
but not vice versa. Any quasitoric manifold admits a canonical smooth
structure (\cite{mmj}).

Submanifolds of the form $\pi^{-1}(F_j), j=1\ldots m$, where $F_j$
is a codimension one facet, are called {\it characteristic
submanifolds}. An {\it omniorientation} of $M^{2n}$ is the
orientation of $M^{2n}$ and all of its characteristic
submanifolds. Every omniorientation determines a canonical
$T^n$-invariant weakly complex structure on $M^{2n}$ (\cite{mmj}),
which may also be constructed in terms of {\it characteristic map}
$\lambda : \m Z^m\to \m Z^n$.

Every smooth projective toric variety admits canonical $T^n$-invariant complex
structure. One may ask the following question (\cite{davis}, Prob. 7.6): find
the criterion for existence of $T^n$-invariant almost
complex structure on quasitoric manifold
$M^{2n}$ in terms of $\lambda$. In this paper we present
the solution of this problem.

$T^n$-invariant weakly complex structure on manifold $M^{2n}$ with
torus action $\alpha$ is an isomorphism $c_{\tau}$ of bundle
$\tau(M^{2n})\oplus \m C^{m-n}$ with complex vector bundle $\xi$.
The map
$$
\xi\stackrel{c_{\tau}^{-1}}{\longrightarrow} \tau(M^{2n})\oplus \m C^{m-n}
\stackrel{d\alpha(t)\oplus I}{\longrightarrow} \tau(M^{2n})\oplus \m C^{m-n}
\stackrel{c_{\tau}}{\longrightarrow}\xi
$$
is complex-linear for every $t$.

If $v$ is a fixed point for $\alpha$, then there is a map $
\tau_v(M^{2n})\stackrel{i}{\rightarrow} \xi_v
\stackrel{\pi'}{\rightarrow} \m C^n, $ where $i$ is an injection and
$\pi'$ is a projection map along fixed subspace $\m C^{m-n}\subset
\xi_v$. The {\it sign} of $v$ is the sign of determinant of the map
$\pi'\circ i$.

An equivalent definition of sign is as follows: we set $\op{sign}(v)=+1$, if
the
orientations of $\tau(M)_v$ determined by orientation of $M$ and orientations of
characteristic submanifolds coincide, and $\op{sign}(v)=-1$, if not.

We call omniorientation {\it positive}, if
signs of all fixed points
are positive.

One can describe notion sign of a fixed point in terms of characteristic
function $\lambda$. Let $e_j$ be a normal vector to facet $F_j$ directed to
interior of $P$. Then $\op{sign}(v)$ is equal to $\det(e_j)\cdot\det\lambda_v$,
where $\lambda_v$ is a
matrix composed of columns of $\lambda$ that correspond to codimension one
facets meeting in $v$.\\

{\bf Theorem 1} {\it Quasitoric manifold $M^{2n}$ admits a $T^n$-invariant almost
complex structure if and only if it admits a positive omniorientation.}\\

Moreover, the following holds.\\

{\bf Theorem 2} {\it Any $T^n$-invariant almost complex structure on $M$ is 
equivariantly equivalent to canonical equivariant weakly complex structure 
\cite{mmj}}.\\

This paper only contains the proof of theorem 1; the proof of theorem 2 will
be published later. (Up to date, the proof exists, but only in Russian).


\section{Remarks and examples}

As follows from \cite{panov}, quasitoric manifold $M$ with positive
omniorientation and corresponding weakly complex structure $J$
satisfies $c_n(J)=\chi(M)$. So by \cite{thomas}, the structure $J$
is equivalent to some almost complex structure on $M$, not
necessarily $T^n$-invariant.

Every smooth projective toric variety $M$ admits canonical
$T^n$-invariant complex structure and symplectic form. So by
\cite{atiyah}, there's a {\it moment map} $\pi:M\to P$, where $P$
is a simple convex polytope. Since all characteristic submanifolds
of $M$ are complex subvarieties, $M$ has positive omniorientation.

Note that in real dimension $6$ there exist examples of
non-projective smooth toric varieties (see \cite{buchpanov}).
Nevertheless, the quotient by
the action of $T^3$ is still isomorphic to simple convex polytope.
Finding an example of non-projective non-quasitoric smooth
toric variety still remains an open problem.

Now we turn to examples of non-algebraic quasitoric manifolds.
Denote by $\m CP^2_k$ a connected sum of $k$ copies of $\m CP^2$
with standart orientation and smooth structure. As follows from
equivariant connected sum constuction (\cite{mmj}), all $\m
CP^2_k$'s are quasitoric. Also, if $k$ is odd, then $\m CP^2_k$ admits a
positive omniorientation. If $k$ is even, then $\m CP^2_k$ doesn't
admit an almost complex structure.\\

{\bf Theorem 2.1 (\cite{wu})} {\it Let $M^4$ be an oriented smooth
manifold with no boundary. Then $M$ admits an almost complex
structure if and only if the number $td(M) =
\frac14(\chi(M)+\op{sign}(M))$ is integer. Here $\op{sign}(M) =
(b_2^+(M) - b_2^-(M))$ is the signature of a cohomology intersection
form on $M$.}\\

We have $td(\m CP^2_k)=\frac{k+1}2$ and this implies non-existence
of almost complex structure on $\m CP^2_k$ for even $k$.

From Theorem 1 we obtain that the manifold $\m CP^2_k$ admits
$T^2$-invariant almost complex structure, if $k$ is odd.
But if $k\geqslant 3$,
$\m CP^2_k$ can't be a toric variety since all toric
varieties are rational and their Todd genus is equal to $1$.

As follows from $\cite{feldman}$, the manifold $\m CP^2_k$ doesn't
admit $T^2$-invariant symplectic structure, if $k\geqslant 2$.\\

{\bf Theorem 2.2 (\cite{feldman})} {\it Let $M$ be a manifold equipped with a
symplectic circle action with only isolated fixed points. Then
$td(M)=1$ if the action is Hamiltonian, and $td(M)=0$ if not.}\\

The proof of following result is based on Seiberg-Witten theory of
 four-dimensional manifolds.\\

{\bf Theorem 2.3 (\cite{taubes})} {\it The manifold $\m CP^2_k$ doesn't
admit a symplectic structure, if $k\geqslant 2$.}\\

Hence, manifolds $\m CP^2_k$ with $k\geqslant 3$, $k$ odd, are
neither symplectic manifolds nor toric varieties -- but they admit
$T^n$-invariant almost complex structure by Theorem 1.

\section{Proof of theorem 1}

The ``only if'' part of Theorem 1 is obvious. If $J$ is an invariant
almost complex structure on $M$, then all characteristic submanifolds
of $M$ are $J$-invariant and all signs of fixed points are positive.
Therefore, it remains to show that existence of positive
omniorientation implies the existence of invariant almost complex
structure.

\subsection{Notations}

Denote by $sk_i(P)$ the set of all $i$-dimensional facets of $P$.
Let $\e : P\to M$ be an embedding of $P$ to $M$ satisfying the
following conditions:
\begin{enumerate}
{\item $\pi\circ\e = id$;}

{\item restriction of $\e$ on $Int\, G$ is smooth for any facet
$G\subset P$.}
\end{enumerate}

An example of such $\e$ is given by composition $P\to P\times T^n\to
M$, where the last map is factorization map (\cite{buchpanov}).

We assume that $M$ admits $T^n$-invariant metric $g$ (see
\cite{bredon}). Then our problem is equivalent to constructing an
operator $J$ on $\tau(M)$ with following properties:
\begin{enumerate}
{\item $J^2=-1$;} {\item $J$ commutes with $d\alpha(t)$ for every $t\in T^n$;}
{\item $J$ is
orthogonal with respect to the metric $g$.}
\end{enumerate}

Recall that the space of orthogonal complex structures on oriented
vector space $\m R^{2i}$ is homeomorphic to $SO(2i)/U(i)$.

The words ``structure on $X$'' will mean ``$T^n$-invariant complex
structure on vector bundle $\tau(M)$ restricted on $X$''.

If $M_G=\pi^{-1}(G)$ is a $2i$-dimensional quasitoric submanifold in
$M$, then we denote by $\xi_1\ldots \xi_{n-i}$ two-dimensional
vector bundles over $M$ corresponding to $(n-i)$ one-dimensional
stationary toric subgroups $T_1\ldots T_{n-i}$ of $M_G$. The fibers
of $\xi_1\ldots \xi_{n-i}$ over $M_G$ are $J$-invariant since $J$
commutes with $d\alpha$.

We'll say that structure $J$ on $\pi^{-1}(Int\,G)$ {\it respects}
omniorientation $o$ if it agrees with
orientations of $\xi_1,\ldots,\xi_{n-i}$ determined by $o$.

Let $V$ be a real oriented Euclidean vector space of even dimension.
Denote by $\m J(V)$ the space of all complex structures of $V$ respecting
orientation and metric. The case of $V=\tau(M_G)|_x$, where $G\subset P$ is a
face, $x\in\e(G)$, will be the most important. Recall that choosing a basis in
$V$ identifies $\m J(V)$ with $SO(2i)/U(i)$, where $i = \dim V$. In particular,
$\m J(V)$ is simply connected.

We'll denote by $\m J_G$ the bundle with fiber $\m J(\tau(M_G)|_x)$ over
$\e(G)$, which is associated with $\tau(M_G)|_{\e(G)}$. Clearly, $\m J_G$ is
trivial. Consider the space $\op{Aut}(\m J(\tau(M_G)|_x))$ of homeomorphisms $\m
J(\tau(M_G)|_x)$ onto itself induced by change of basis in the space
$\tau(M_G)|_x$.

Let us fix an arbitrary trivialization of $\tau(M_G)$ over $\e(G)$. Then it
determines also a trivialization of associated bundle $\m J_G$. Any other
trivialization of $\m J_G$ is then determined by an arbitrary continious map
$G\to Aut(\m J(\tau(M_G)|_x))$, where $x\in\e(G)$ is some fixed point. Since
$SO(2i)$ is connected, $\op{Aut}(\m J(\tau(M_G)|_x))$ is also connected.
We obtain that the space of trivializations of $\m J_G$ determined by
trivialization of $\tau(M_G)$ over $\e(G)$ is also connected.

\subsection{Positivity and one-dimensional facets}

\bigskip

{\bf Lemma 3.1} {\it If $M_G\subset M$ is a quasitoric submanifold, then
$\xi_j\,\bot\, M_G$ and $\xi_j\,\bot\,\xi_k$, if $j\ne k$.}\\

$\Box$ Let $v\in \tau(M_G), v_j\in\xi_j$ and $v_k\in \xi_k$ be
nonzero vectors in the point $x\in M_G$, and $t_{\pi}\in T_j$ be an
element of toric subgroup $T_j$ corresponding to multiplication by
$-1$. Since $T_j$ acts trivially on $\xi_k$ and $\tau(M_G)$, we have
$g(v,v_j) = g(v,t_{\pi}v_j) = g(v,-v_j)=0$ and $g(v_j,v_k) =
g(t_{\pi}v_j,v_k) = g(-v_j, v_k) =0$. $\Box$

Let $o$ be any omniorientation of $M$, not necessarily positive.
Then $o$ determines structure $J$ on $\pi^{-1}(sk_0(P))$ in the
following way: $J$ is rotation by angle $\pi/2$ in the direction
specified by $o$ on every fiber $\xi_1|_x,\ldots,\xi_n|_x$, where $x\in
\pi^{-1}(sk_0(P))$.\\

{\bf Proposition 3.2} {\it Let $J$ be a structure on $\pi^{-1}(sk_0(P))$
determined by $o$. Then $J$ may be extended over $\pi^{-1}(sk_1(P))$
if and only if $o$ is positive.}\\

$\Box$ Let us prove the ``if'' part first. Consider $I\subset P$ --
an edge connecting vertices $x_0$ and $x_1$, $T_1\ldots T_{n-1}$ --
toric subgroups corresponding to $I$. Then $\tau(M)\simeq
\tau(\pi^{-1}(I))\oplus\xi_1\oplus\ldots\oplus\xi_{n-1}$ over
$\pi^{-1}(I)$. The orientation of bundles $\xi_1\ldots\xi_{n-1}$
over $\pi^{-1}(I)$ is determined by $o$.

Denote by $W_0$ and $W_1$ orthogonal complements to
$\xi_1\ldots\xi_{n-1}$ at the points $\pi^{-1}(x_0)$ and $\pi^{-1}(x_1)$
respectively. Since $J$ was constructed by $o$, $J$ is a rotation by $\pi/2$ in
spaces
$W_0$ and $W_1$. Note that $W_0$ and $W_1$ are tangent spaces to $\pi^{-1}(I)$
at $\pi^{-1}(x_0)$ and $\pi^{-1}(x_1)$. The positivity of $o$ implies that
orientations of $W_0$ and $W_1$ agree as orientations of
$\tau(\pi^{-1}(I))=\tau(S^2)$. So we can define $J$ on $\tau(\pi^{-1}(I))$
as rotation by $\pi/2$ in the direction specified by orientations of $W_0$ and
$W_1$. On the bundles $\xi_1,\ldots\xi_{n-1}$ structure $J$ is defined as
rotation by $\pi/2$ in the direction specified by $o$.

The "only if" part of Prop.3.2 is proved in similar fashion.
If $J$ is defined on entire $\pi^{-1}(I)$, then the orientations of $M$ at
$\pi^{-1}(x_0)$ and $\pi^{-1}(x_1)$ agree -- and that holds for every pair of
adjacent vertices $(x_0,x_1)$. $\Box$

\subsection{Triviality of higher obstructions}

Now we consider the case $i>1$. Suppose that $J$ is defined on
$\pi^{-1}(sk_{i-1}(P))$ and we're trying to extend it over
$\pi^{-1}(sk_{i}(P))$.\\

{\bf Lemma 3.3} {\it Let $G\subset P$ be an $i$-dimensional facet. The space of
positive structures $J$ over $\pi^{-1}(Int\,G)$ is homeomorphic to space of
continious mappings $\op{Map}(Int\,G, \m J(\tau(M_G)|_x))$, where
$x\in\e(Int\,G)$ is a fixed point.}\\

$\Box$ The positivity of $J$ implies that $J$ is uniquely defined on bundle
$\tau(\pi^{-1}(Int\,G))^{\bot}\simeq \xi_1\oplus\ldots\oplus\xi_{n-i}$.
So it suffices to define $J$ on the tangent
vector bundle $\tau(\pi^{-1}(Int\,G))\simeq \tau(M_G)|_{Int\,G}$.

Denote by $T^i$ any $i$-dimensional toric subgroup complementary to $T^{n-i}$ in
$T^n$. Then $T^i$ acts freely on $\pi^{-1}(Int\,G)$ and the quotient by the
action is homeomorphic to $Int\,G$ itself. Since $J$ is $T^n$-invariant
structure, it suffices to define $J$ only on $\e(Int\,G)$. The bundle
$\tau(M_G)$ is trivial oriented $2i$-dimensional real vector bundle
over $\e(Int\, G)$ -- and so $J$ may be chosen in an arbitrary way over it (the
only condition is that $J$ must be orthogonal with respect to the metric $g$).
If we fix some trivialization of associated bundle $\m J_G$ over $\e(Int\,G)$,
then $J$ is determined by an arbitrary continious map $Int\,G\to\m
J(\tau(M_G)|_x))$, where $x\in\e(G)$ is a fixed point. $\Box$

Note that the bundle $\tau(M_G)$ over $\e(G)$ is also trivial. This implies
existence of a canonical isomorphism of homotopy groups $\pi_*(\m
J(\tau(M_G)|_x))$ and $\pi_*(\m
J(\tau(M_G)|_y))$ for any points $x,y\in\e(G)$. Recall that the homotopy groups $\pi_*(\m
J(\tau(M_G)|_x))$ are independent of choice of starting point if $\m
J(\tau(M_G)|_x))$, since $J(\tau(M_G)|_x))$ is simply connected.

Let us fix a point $x\in\e(G)$ and trivialization of bundle
$\tau(M_G)$ over $\e(G)$. Then the bundle $\m J_G$ is also trivialized. Since
$J$ is already defined over $\tau(M)|_{\e(\partial G)}$ and
$\tau(M_G)\subset\tau(M)$ is $J$-invariant subbundle, we obtain a continious map
$f:\partial G\to \m J(\tau(M_G)|_x)$. Denote by $C_G$ the homotopy class of
spheroid $f$ in $\pi_{i-1}(J(\tau(M_G)|_x))$.

Now we'll show that the class $C_G$ is well-defined. Clearly, $C_G$ is
independent of the choice of $x$. Suppose that we have changed
trivialization of $\tau(M_G)$ over $\e(G)$. Then new trivialization is given by
a continious map
$\phi:G\to \op{Aut}(\m J(\tau(M_G)|_x))$ (see subsection 3.1). The
map $f$ is changed to $g: \partial G\to J(\tau(M_G)|_x))$, where $g(y) =
\phi(y)\circ f(y)$.

Let $\phi_t:G\to \op{Aut}(\m J(\tau(M_G)|_x))$ be a continious family of maps
satisfying $\phi_0=id,\,\phi_1=\phi$ (recall that
$\op{Aut}(\m J(\tau(M_G)|_x))$ is connected). Then family of maps
$\phi_t(y)\circ f(y)$ provides a homotopy between $f$ and $g$. This completes
the proof of that $C_G$ is well-defined.

Note that $C_G$ depends on the orientation of $G$.

We summarize all in the following statement.\\

{\bf Lemma 3.4} {\it Let $J$ be a structure on
$\pi^{-1}(sk_{i-1}(P))$ respecting $o$. Then $J$ may be extended over
$\pi^{-1}(sk_{i-1}(P))\cup M_G$ if and only if $C_G=0$ in the
group $\pi_{i-1}(SO(2i)/U(i))$.}\\

Every polytope $P$ has a canonical cellular decomposition, with facets playing
the role of cells.
We define cellular cochain $\sigma^i_J\in
C^i(P,\pi_{i-1}(SO(2i)/U(i))$ by the rule $\sigma^i_J(G)=C_G$. By definition,
$\sigma^i_J$ is zero if and only if $J$ may be extended from
$\pi^{-1}(sk_{i-1}(P))$ to $\pi^{-1}(sk_i(P))$.

Strictly speaking, to call $\sigma^i_J$ a cocycle we have to identify all
homotopy groups $\pi_{i-1}(\m J(\tau(M_G)|_x))$ for all $i$-dimensional facets
$G\subset P$.\\

{\bf Lemma 3.5} {\it Let $\dim G=i$, $j\leqslant 2i-2$, $x\in\e(G)$, $y\in\e(P)$.
Then homotopy groups $\pi_j(\m J(\tau(M_G)|_x))$ and $\pi_j(\m
J(\tau(M)|_y))$ are canonically isomorphic.}\\

$\Box$ Note that lemma 3.5 implies that $\sigma^i_J$ is a cochain, since
$i-1\leqslant 2i-2$. Since homotopy groups $\pi_j(\m J(\tau(M)|_x))$ and
$\pi_j(\m J(\tau(M)|_y))$ are canonically isomorphic, one can assume that $x=y$.

Consider an arbitrary embedding of two facets $H\subset L$ such that $\dim
H\geqslant i,\, \dim L=\dim H+1$. Let $x\in\e(H)$. Then there is an embedding of
spaces of complex structures
$$
C(H,L):\m J(\tau(M_H)|_x)\to \m J(\tau(M_L)|_x),
$$
defined by formula $J\to J\oplus t_{\pi/2}$. Here $t_{\pi/2}$ is a rotation by
$\pi/2$ in two-dimensional orthogonal complement $\tau(M_H)^{\bot}\subset
\tau(M_L)$. The direction of $t_{\pi/2}$ is specified by coorientation of
$\tau(M_H)|_x$ in $\tau(M_L)|_x$.

IF we fix the basis in $\tau(M_H)|_x$ and then add two vectors
to obtain a basis in $\tau(M_L)|_x$, then $c(H,L)$ turns into a canonical
embedding of homogeneous spaces $SO(2r)/U(r)\to SO(2r+2)/U(r+1)$, where $r=\dim
H$.\\

{\bf Lemma 3.6} {\it The map
$$
c_*:\pi_{j}(SO(2r)/U(r))\to \pi_{j}(SO(2r+2)/U(r+1))
$$

is an isomorphism for every $j\leqslant 2r-2$.}\\

$\Box$ The embedding $c$ may be viewed as the composition
$$
SO(2i)/U(i)\to SO(2i+2)/U(i)\to SO(2i+2)/U(i+1),
$$
where first map is an embedding of fiber to bundle space over
$SO(2i+2)/SO(2i)$ and the second map is a projection map of bundle
space with fiber $S^{2i+1}$. An exact homotopy sequence now implies
that $c$ induces isomorphism of corresponding homotopy groups up to
dimension $2i-2$. $\Box$

Now we consider an arbitrary chain of embeddings $G =
G_0\subset\ldots\subset G_{n-i}=P$, where $\dim G_{s+1} = \dim
G_s+1$ for every $s$. Define an isomorphism $c_*(G,P): \pi_j(\m
J(\tau(M_G)|_x))\to \pi_j(\m J(\tau(M)|_x))$ by formula $c_*(G,P) =
c_*(G_{n-i-1},G_{n-i})\circ\ldots\circ ё_*(G_0, G_1)$. Then
$c_*(G,P)$ is just what we need; it suffices to prove following
statement.\\

{\bf Lemma 3.7} {\it Isomorphism $c_*(G,P)$ is independent on the chain
$G_0\subset\ldots\subset G_{n-i}$.}\\

$\Box$ Consider an arbitrary subchain of the form $G_{s-1}\subset
G_s\subset G_{s+1}$. Then there exists a unique facet $Q$ such that
$G_{s-1}\subset Q\subset G_{s+1}$ and $G_s\ne Q$. If $x\in G_{s-1}$,
then the following diagram of embeddings
$$
\begin{CD}
\m J(\tau(M_{G_{s-1}})|_x)) @>>> \m J(\tau(M_{G_s})|_x))\\
@VVV @VVV\\
\m J(\tau(M_Q)|_x)) @>>> \m J(\tau(M_{G_{s+1}})|_x)
\end{CD}
$$
is commutative. So if we replace $G_s$ with $Q$, the resulting
isomorphism $c_*(G,P)$ won't change.

We'll call two chains connecting $G$ and $P$ {\it equivalent}, if
one may be obtained from another by sequence of operations we've
just described. Let us prove by induction on $n-i$ that any two
chains connecting $G$ and $P$ are equivalent. The base statement
($n-i=1$) is obvious since the chain is unique. Now consider two
different chains $G=G_0^1\subset\ldots\subset G_{n-i}^1=P$ and
$G=G_0^2\subset\ldots\subset G_{n-i}^2=P$. Let $Q\subset P$ be a
facet satisfying $\dim Q=i+2$, $G_1^1\subset Q$ and $G_1^2\subset
Q$. Consider an arbitrary chain $\zeta$ connecting $Q$ and $P$. By
induction hypothesis, two chains $G_1^1\subset
G_2^1\subset\ldots\subset P$ and $G_1^1\subset Q\subset\ldots\subset
P$, where $Q$ is connected with $P$ by $\zeta$, are equivalent. This
means $G\subset G_1^1\subset G_2^1\subset\ldots\subset P$ is
equivalent to the chain $G\subset G_1^1\subset Q\subset\ldots\subset
P$. Moreover, $G\subset G_1^1\subset Q\subset\ldots\subset P$ is
equivalent to $G\subset G_2^1\subset Q\subset\ldots\subset P$, due
to the choice of $Q$. Finally, applying the induction hypothesis
once more, we obtain that $G\subset G^2_1\subset
Q\subset\ldots\subset P$ and $G\subset G^2_1\subset
G_2^2\subset\ldots\subset P$ are equivalent. $\Box$

This lemma completes the proof of that $\sigma^i_J$ is well-defined.
We summarize this in the following statement.\\

{\bf Lemma 3.8} {\it Let $J$ be a structure on $\pi^{-1}(sk_{i-1})$
respecting $o$. Then one can define an obstruction cochain
$\sigma^i_J\in C^i(P,\pi_{i-1}(SO(2i)/U(i))$ which is zero iff $J$
may be extended to $\pi^{-1}(sk_i(P))$.}\\

Our next aim is to prove that $\sigma^i_J$ is a cocycle.\\

{\bf Lemma 3.9.} {\it Suppose that $J$ is a structure that respects $o$
and is defined on $\pi^{-1}(sk_{i-1}(P))$, $Q$ is an
$(i+1)$-dimensional facet of $P$. Then
$$
\sum\limits_{G\subset \partial Q} \sigma^i_J(G) = 0.
$$}

$\Box$ Lemma 3.5 guarantees that we can use a single notation
$c_*(G,Q)$ for all isomorphisms $c_*(G,Q):\pi_{i-1}(\m
J(\tau(M_G)|_x))\to\pi_{i-1}(\m J(\tau(M_Q)|_y))$, $x\in
\e(G),\,y\in \e(Q)$.

The structure $J$ on the bundle $\tau(M_G)|_{\e(G)}$ automatically
defines structure $c(J)$ on $\tau(M_Q)|_{\e(G)}$ by the formula
$J\to J\oplus t_{\pi/2}$. The bundle $\tau(M_Q)$ may be trivialized
over $\e(Q)$. Since $\e(G)\subset\e(Q)$ for every $G$, there is a
well-defined continious map $f_J:Q\cap sk_{i-1}(P)\to \m
J(\tau(M_Q)|_y)$, $y\in\e(Q)$ a fixed point. By its construction,
homotopy class of $f_J|_{\partial G}$ is equal to $c_*\sigma^i_J$.

Let $U_{\eps}$ be closed $i$-dimensional tubilar neighbourhood of
$Q\cap sk_{i-1}(P)$ in $\partial Q$ (such that $G\cap
U_{\eps}\ne\emptyset$ for every $G$). The set $Q\cap sk_{i-1}(P)$ is
a deformational retract of $U_{\eps}$, and the retraction determines
a map $\tilde f_J: U_{\eps}\to SO(2i+2)/U(i+1)$. The space
$U_{\eps}$ is homeomorphic to $i$-dimensional sphere with $k$ holes;
let $S_1,\ldots,S_k$ be the boundary spheres of the these holes.

$$
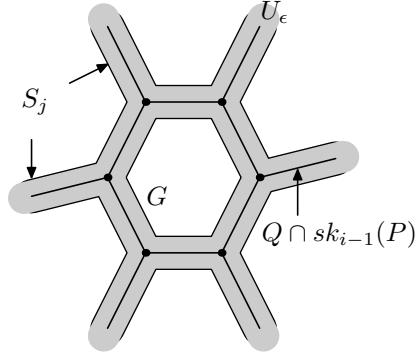
\epsfbox{am1.13}$$

Then any homotopy class $\tilde f_J|_{S_j}$ coincides with
$c_*\sigma^i_J(G)$, if $S_j\subset G$. Since the sum of spheroids
$\tilde f_J|_{S_j}$ must be zero in $\pi_{i-1}(SO(2i+2)/U(i+1))$,
this completes the proof. $\Box$\\

{\bf Lemma 3.10} {\it Suppose that $J$ is a structure on
$\pi^{-1}(sk_{i-1}(P))$ that respects $o$ and $\sigma^i_J$ is a
coboundary. Then one can change $J$ on
$\pi^{-1}(sk_{i-1}(P)\setminus sk_{i-2}(P))$ and obtain new
structure $J'$ on $\pi^{-1}(sk_{i-1}(P))$ such that $\sigma^i_{J'}=0$.}\\

$\Box$ By Lemma 3.6 the map $c:SO(2i-2)/U(i-1)\to SO(2i)/U(i)$
induces an isomorphism of homotopy groups up to dimension $2i-4$. If
$i>2$, then $i-1\leqslant 2i-4$, and if $i=2$, then $SO(4)/U(2)$ and
$SO(2)/U(1)$ are simply connected -- so the map $c_*:
\pi_{i-1}(SO(2i-2)/U(i-1))\to \pi_{i-1}(SO(2i)/U(i))$ is an
isomorphism for $i>1$.

Let $\sigma^i_J=\delta\beta$ and $H$ be some $(i-1)$-dimensional
facet. The proof of lemma 3.5 implies that we can use a single
notation $c_*$ for every isomorphism of the form
$c_*(H,G_i):\pi_{i-1}(\m J(\tau(M_H)|_x))\to\pi_{i-1}(\m
J(\tau(M_{G_i})|_y))$, where $x\in H$, $y\in G_i$, $\dim G_i = i$.
Now apply lemma 3.3: the space of structures on $\pi^{-1}(Int\,H)$
respecting $o$ is homeomorphic to space of continious maps $Int\,H
\to \m J(\tau(M_H)|_x)$. Let $f_H: Int\,H \to \m J(\tau(M_H)|_x)$ be
the map corresponding to $J$.

Let us identify $Int\,H$ with open $(i-1)$-dimensional ball in $\m
R^i$ of radius $1$. Consider the map $\tilde f_H: Int\,H \to \m
J(\tau(M_H)|_x)$ satisfying following conditions:
\begin{itemize}
{\item $\tilde f_H(x) = f_H((2|x|-1)\cdot x)$, if $1/2\leqslant
|x|\leqslant 1$;} {\item the homotopy class of spheroid $\tilde
f_H|_{\{|x|\leqslant 1/2\}}$ is $(-c_*^{-1}\beta(H))$:}
\end{itemize}

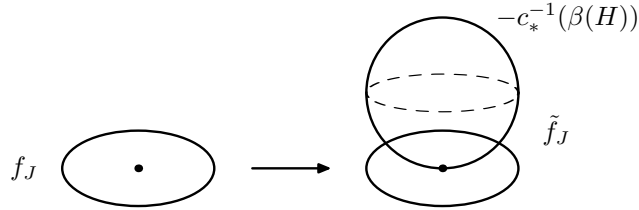
$$ \epsfbox{am1.12} $$

If we replace $f_H$ with $\tilde f_H$, we'll obtain new structure
$J'$ on $\pi^{-1}(sk_i(P))$. Then $\sigma^i_{J'}(G) = \sigma^i_J(G)
- \beta(H)$ for every $i$-dimensional facet $G$ such that $H\subset
G$. So if we replace $f_H$ with $\tilde f_H$ for every
$(i-1)$-dimensional facet $H\subset P$, we'll have
$\sigma^i_{J'}=0$, since $\sigma^i_J=\delta\beta$. $\Box$

\subsection{Proof of the main theorem}

Now we can prove theorem 1. Recall that we construct $J$ by
induction on $i$, starting from $\pi^{-1}(sk_0(P))$. Proposition 3.2
and the positivity of $o$ guarantee that $J$ may be extended from
$\pi^{-1}(sk_0(P))$ to $\pi^{-1}(sk_1(P))$. Lemmas 3.8 and 3.9 imply
that for $i>1$ obstruction to extending $J$ from
$\pi^{-1}(sk_{i-1}(P))$ to $\pi^{-1}(sk_i(P))$ is a cocycle
$\sigma^i_J\in C^i(P,\pi_{i-1}(SO(2i)/U(i)))$. Since $P$ has trivial
homology, $\sigma^i_J$ is a coboundary, so by Lemma 3.10 there
exists structure $J'$ on $\pi^{-1}(sk_{i-1}(P))$ that may be
extended to $\pi^{-1}(sk_i(P))$.

\section{Acknowledgements}

The author is grateful to his advisor, Prof.V.M.Buchstaber, and to
T.E.Panov for constant attention to this work, and also to
N.E.Dobrinskaya for useful e-mail correspondence. The statement of
Theorem 1 was proved earlier by N.E.Dobrinskaya for $n\leqslant 7$.


\begin{thebibliography}{9999}

\bibitem{wu} W.-T. Wu. Sur le classes caracteristique des structures fibrees
spheriques // Actualites Sci. Industr. 1183 (1952).
\bibitem{thomas} E. Thomas. Complex Structures on Real Vector Bundles //
Amer. J. Math. 1967. V.89. P.887-908.
\bibitem{bredon} Glen E. Bredon. Introduction to Compact Transformation Groups.
Pure and Applied Math. 46, Academic Press, 1972.
\bibitem{atiyah} M. F. Atiyah. Convexity and Commuting Hamiltonians // Bull.
London Math. Soc. 1982. V.14, N1. P. 1-15.
\bibitem{gs} V. Guilemin, S. Sternberg. Convexity properties of the moment
mapping // Inv. Math. 67 (1982), N3, P. 491-513.
\bibitem{davis} M. Davis, T. Januskiewicz. Convex Polytopes, Coxeter Orbifolds
and Torus Actions // Duke Math J. 1991. V.62 N2 P. 417-451.
\bibitem{taubes} C.H.Taubes. The Seiberg-Witten Invariants and Symplectic Forms
// Math. Res. Lett. 1994. V.1. P.809-822.
\bibitem{panov} T. E. Panov. Hirzebruch Genera of Manifolds with Torus Action //
Izv. RAN. Ser. Mat., 65:3 (2001), 123Ц138
\bibitem{feldman} K.E.Feldman. Hirzebruch Genera of Manifolds Equipped with a
Hamiltonian Circle Action. arXiv:math/0110028v2
\bibitem{buchpanov} Victor M. Buchstaber, Taras E. Panov.
Torus Actions and Their Applications in Topology and Combinatorics.
AMS Bookstore, 2002.
\bibitem{mmj} V. Buchstaber, T. Panov, N. Ray. Spaces of Polytopes and Cobordism
of Quasitoric Manifolds // Moscow Math. J (2007) V.7 N2.



\end{thebibliography}
\end{document}